\documentclass{article}
\catcode`\@=11
\newtheorem{theorem}{Theorem}
\newtheorem{lemma}{Lemma}
\newtheorem{proposition}{Proposition}

\newtheorem{example}{Example}
\newtheorem{remark}{Remark}

\newtheorem{note}{Note}
\newtheorem{corollary}{Corollary}
\def\demo{\noindent{\bf Proof .-}}
\def\section{\@startsection {section}{1}{\z@}{-3.5ex plus -1ex
minus-.2ex}{2.3ex plus .2ex}{\normalsize\bf}}

\begin{document}
\begin{center}
{\Large\bf \textsc{On the arithmetical rank of an intersection of ideals}}
\end{center}
\vskip.5truecm
\begin{center}
{Margherita Barile\footnote{Partially supported by the Italian Ministry of Education, University and Research.}\\ Dipartimento di Matematica, Universit\`{a} di Bari, Via E. Orabona 4,\\70125 Bari, Italy}\footnote{e-mail: barile@dm.uniba.it, Fax: 0039 080 596 3612}
\end{center}
\vskip1truecm
\noindent
{\bf Abstract} We determine sets of elements which, under certain conditions, generate an intersection of ideals up to radical.   

\section*{Introduction} Let $R$ be Noetherian commutative ring with identity. We say that some elements  $\gamma_1,\dots, \gamma_r\in R$ generate an ideal $I$ of $R$  {\it up to radical} if $\sqrt{I}=\sqrt{(\gamma_1,\dots, \gamma_r)}$. The smallest $r$ with this property is called the {\it arithmetical rank} of $I$, denoted ara\,$I$. It is well known that height\,$I\leq$\,ara\,$I$. If equality holds, $I$ is called a {\it set-theoretic complete intersection} ({\it s.t.c.i.}). Determining the arithmetical rank of an ideal $I$, or at least a satisfactory upper bound for it, is, in general, a hard task. In this paper we give some criteria for the case where $I$ is the intersection of two ideals, whose generators are linked by special divisibility conditions. Our results yield constructive methods which can be applied to certain polynomial ideals generated by squarefree monomials. In particular, we characterize a class of monomial ideals of minimal multiplicity which are s.t.c.i..\newline 
Other results on the arithmetical rank of intersections of ideals can be found in \cite{B4}. \newline

\section{The Main Theorem}
For the proof of our main theorem we will need the following preliminary result, which was presented in \cite{B1}, and is a generalization of the lemma in \cite{SV}, p.~249. 
\begin{lemma}\label{lemma1} Let $P_1,\dots, P_r$ be finite subsets of $R$, and set $P=\bigcup_{i=1}^rP_i$. Suppose that
\begin{list}{}{}
\item[(i)] $P_1$ has exactly one element;
\item[(ii)] if $p$ and $p'$ are different elements of $P_i$ $(1<i\leq r)$ then $(pp')^m\in\left(\bigcup_{j=1}^{i-1}P_j\right)$ for some positive integer $m$. 
\end{list}
\noindent
We set $q_i=\sum_{p\in P_i}p^{e(p)}$, where $e(p)\geq1$ are arbitrary integers. Then we get
$$\sqrt{(P)}=\sqrt{(q_1,\dots,q_r)}.$$
\end{lemma}
Before stating our main result, we introduce the following notation: given $\alpha_1,\dots, \alpha_t$, $\beta_1,\dots, \beta_t\in R$, we write $[\alpha_1,\dots, \alpha_t]\subset [\beta_1,\dots, \beta_t]$ if, for all indices $i$ such that $1\leq i\leq t$, $(\alpha_1,\dots, \alpha_i)\subset (\beta_1,\dots, \beta_i)$.\par\smallskip\noindent

\begin{theorem}\label{main} Let $\mu_1,\dots, \mu_u, \nu_1,\dots, \nu_v$ be elements of $R$ and consider the ideals $I=(\mu_1,\dots, \mu_u)$, and $J=(\nu_1,\dots, \nu_v)$ of $R$. Suppose that there is a positive integer $s$ with $s<u$ and $s\leq v$ such that
\begin{equation}\label{0}[\mu_1,\dots, \mu_s]\subset[\nu_1,\dots, \nu_s].\end{equation}
\noindent
 For  $1\leq i\leq u-s$ set
\begin{equation}\label{1}\gamma_i=\sum_{j=s+1}^{s+i}\mu_j\nu_{i-j+s+1},\end{equation}
\noindent
and, for $u-s+1\leq i\leq u$, set
\begin{equation}\label{2}\gamma_i=\mu_{i-u+s}+\sum_{j=s+1}^{s+i}\mu_j\nu_{i-j+s+1}.\end{equation}
\noindent
If $s\leq v-2$, for $u+1\leq i\leq u+v-s-1$, set further
\begin{equation}\label{3}\gamma_i=\sum_{j=s+1}^{s+i}\mu_j\nu_{i-j+s+1}.\end{equation}
\noindent
In equations (\ref{1}), (\ref{2}) and (\ref{3}) we have adopted the following convention: we have set $\mu_h=0$ whenever $h>u$ and $\nu_k=0$ whenever $k>v$. \newline
Then, if $s\leq v-2$, 
\begin{equation}\label{4}\sqrt{I\cap J}=\sqrt{(\gamma_1,\dots, \gamma_{u+v-s-1})},\end{equation}
\noindent
and, if $v-1\leq s\leq v$, 
\begin{equation}\label{5}\sqrt{I\cap J}=\sqrt{(\gamma_1,\dots, \gamma_u)}.\end{equation}
\end{theorem}
\begin{note} {\rm We can rewrite (\ref{1}) and (\ref{3}) in the form:
$$\gamma_i=\sum_{{h+k=i+s+1}\scriptstyle\atop{h\geq s+1}}\mu_h\nu_k,\qquad\qquad(i\leq u-s\mbox{ or } i\geq u+1),$$
\noindent
and (\ref{2}) in the form
$$\gamma_i=\mu_{i-u+s}+\sum_{{h+k=i+s+1}\scriptstyle\atop{h\geq s+1}}\mu_h\nu_k.\qquad\qquad(u-s+1\leq i\leq u)$$
\noindent
We are going to use this simplified notation in the proof of Theorem \ref{main}.
}\end{note}
\demo We prove the equalities (\ref{4}) and (\ref{5}) by showing the two inclusions. For all indices $i$, all summands of $\gamma_i$ belong to $\sqrt{I\cap J}$: this is certainly true for the terms of the form $\mu_h\nu_k$; moreover, for $u-s+1\leq i\leq u$, we have that $1\leq i-u+s\leq s$, so that, by (\ref{0}), $\mu_{i-u+s}\in I\cap J$. Thus inclusion $\supset$ is true. We now prove $\subset$. For all indices $i$, let $P_i$ be the set of summands of $\gamma_i$. Here we are excluding the terms containing $\mu_h$ with $h>u$ or $\nu_k$ with $k>v$. We show that the sets $P_i$ fulfill the assumption of Lemma \ref{lemma1}. Since $P_1=\{\mu_{s+1}\nu_1\}$, (i) is satisfied. Now let $i\geq 2$ and let $p$ and $p'$ be two distinct elements of $P_i$. Then one of the following two cases occurs:
\begin{list}{}{}
\item[(a)] $p=\mu_h\nu_k$ and $p'=\mu_{h'}\nu_{k'}$, where, without loss of generality, we may assume that $1\leq k'< k$. Let $l=k-k'$. Note that 
$pp'$ is divisible by $p''=\mu_h\nu_{k'}$. Since $h+k=i+s+1$, and $h\geq s+1$,  we have $k\leq i$, so that $l<i$. Moreover, $h+k'=h+k-l=i-l+s+1$, so that $p''\in P_{i-l}$, and, consequently, $pp'\in P_{i-l}$. 
\item[(b)] $u-s+1\leq i\leq u$ and $p=\mu_{i-u+s}$, $p'=\mu_{h'}\nu_{k'}$. Then $1\leq i-u+s\leq s$, so that, by (\ref{0}), $pp'\in(\nu_1,\dots, \nu_{i-u+s})(\mu_{h'}\nu_{k'})$. Moreover, for $1\leq k\leq i-u+s$, $\mu_{h'}\nu_k\in P_{h'+k-s-1}$, where $h'+k-s-1\leq h'+i-u-1\leq u+i-u-1=i-1$. Hence $pp'\in \left(\cup_{j=1}^{i-1}P_j\right)$. 
\end{list}
This shows that (ii) is fulfilled, too. It follows that
$$\sqrt{(P)}=\cases{\sqrt{(\gamma_1,\dots, \gamma_{u+v-s-1})}&if $s\leq v-2$,\cr\cr
\sqrt{(\gamma_1,\dots, \gamma_u)}&if $v-1\leq s\leq v$.}
$$
\noindent
Now it remains to show that 
\begin{equation}\label{7} \sqrt{I\cap J}\subset\sqrt{(P)}.\end{equation}
\noindent
Recall that $\sqrt{I\cap J}=\sqrt{IJ}$, and  $IJ$ is generated by all products $\mu_h\nu_k$, with $1\leq h\leq u$, $1\leq k\leq v$. If $h\geq s+1$, then $\mu_h\nu_k\in P_{h+k-s-1}$: note that $h+k-s-1\leq u+v-s-1$. If $1\leq h\leq s$, then $\mu_h\in P_{h+u-s}$. In either case, $\mu_h\nu_k\in (P)$. This shows (\ref{7}) and completes the proof. 
\par\bigskip\noindent
As an immediate consequence of Theorem \ref{main} we have:
\begin{corollary}\label{corollary} Let $\mu_1,\dots, \mu_u, \nu_1,\dots, \nu_v$ be elements of $R$ and consider the ideals $I=(\mu_1,\dots, \mu_u)$, and $J=(\nu_1,\dots, \nu_v)$ of $R$. Suppose that there is a positive integer $s$ with $s<u$ and $s\leq v$ such that
$[\mu_1,\dots, \mu_s]\subset[\nu_1,\dots, \nu_s]$. Then 
$$r={\rm ara}\,(I\cap J)\leq \cases{u+v-s-1&if $s\leq v-2$,\cr
u&if $v-1\leq s\leq v$.}
$$
\noindent
Moreover, the elements $\gamma_1,\dots, \gamma_r\in R$ generating $I\cap J$ up to radical can be chosen in such a way that
\begin{equation}\label{statement1}[\gamma_1,\dots, \gamma_u]\subset[\mu_{s+1},\dots, \mu_{u}, \mu_1,\dots, \mu_s],\end{equation}
and\noindent
\begin{equation}\label{statement2}[\gamma_1,\dots, \gamma_v]\subset[\nu_1,\dots, \nu_v].\end{equation}
\end{corollary}
\demo We only prove statement (\ref{statement2}). For all indices $i$ such that $1\leq i \leq u-s$ from (\ref{1}) we deduce that
$\gamma_i\in(\nu_1\,\dots, \nu_i).$
Now let $u-s+1\leq i\leq u$. Then, by (\ref{2}), 
$$\gamma_i\in(\mu_{i-u+s},\nu_1,\dots, \nu_i ),$$
\noindent
where $1\leq i-u+s\leq s$, and, moreover, $i-u+s<i$. Hence, in view of (\ref{0}), 
$\mu_{i-u+s}\in(\nu_1,\dots, \nu_i),$
 so that $\gamma_i\in(\nu_1, \dots, \nu_i).$ If $v\leq u$, this yields the claim. Otherwise, since $v\leq u+v-s-1$, for all indices $i$ with $u+1\leq i \leq v$, by (\ref{3}) we have that $\gamma_i\in(\nu_1,\dots, \nu_i)$. This completes the proof of (\ref{statement2}). Statement (\ref{statement1}) can be shown in a similar way.
\begin{remark}\label{remark1}{\rm The claim of Corollary \ref{corollary} also holds for $s=0$:  it is well known that, whenever $I$ and $J$ are ideals of $R$ generated by $u$ and $v$ elements respectively, then 
$$r={\rm ara}\,(I\cap J)\leq u+v-1.$$
\noindent
A proof of this result and its generalization to the intersection of any finite number of ideals is given in \cite{SV}, Theorem 1.}
\end{remark}
\begin{example}\label{example1}{\rm Let $K$ be a field. Consider the following ideal of $R=K[x_1,\dots, x_5]$:
$$I=(x_1, x_3, x_5)\cap(x_2, x_4, x_5)\cap (x_1, x_4, x_5)\cap(x_2, x_3, x_5)\cap(x_1, x_3, x_6)\cap(x_2, x_5, x_7),$$
\noindent
which is of pure height 3. We have  $I=(x_1x_2, x_1x_5, x_3x_5, x_5x_6, x_3x_4x_7, x_2x_3x_4)$. We show that ara\,$I=3$, i.e., $I$ is a s.t.c.i.. Note that $I=I_1\cap I_2\cap I_3$, where 
$$I_1=(x_1x_2, x_3x_4, x_5),\ I_2=(x_1, x_3, x_6),\ I_3=(x_2, x_5, x_7).$$
\noindent First set
$$\begin{array}{lll} 
\mu_1=x_1x_2&& \nu_1=x_1\\
\mu_2=x_3x_4&& \nu_2=x_3\\
\mu_3=x_5&& \nu_3=x_6.
\end{array}
$$
\noindent
Then Corollary \ref{corollary} applies to $I_1\cap I_2$ with  $u=v=3$, $s=2$ and $r=3$. Hence $I_1\cap I_2$ is generated up to radical by the following 3 elements:
$$\gamma_1=x_1x_5,\ \gamma_2=x_1x_2+x_3x_5,\ \gamma_3=x_3x_4+x_5x_6.$$
\noindent
Now set 
$$\begin{array}{lll} 
\mu_1'=\gamma_1&& \nu'_1=x_5\\
\mu_2'=\gamma_2&& \nu'_2=x_2\\
\mu_3'=\gamma_3&& \nu'_3=x_7.
\end{array}
$$
\noindent
Then Corollary \ref{corollary} applies to the intersection of the ideals $(\gamma_1, \gamma_2, \gamma_3)$ and $I_3$ with $u=v=3$ and  $s=2$. Thus this intersection is generated up to radical by 
\begin{eqnarray*}
\gamma'_1&=&(x_3x_4+x_5x_6)x_5,\\
\gamma'_2&=&x_1x_5+(x_3x_4+x_5x_6)x_2,\\
 \gamma'_3&=&x_1x_2+x_3x_5+(x_3x_4+x_5x_6)x_7.
\end{eqnarray*}
\noindent
These elements generate $I$ up to radical.}
\end{example}
\begin{example}\label{example2} {\rm In the polynomial ring $R=K[x_1,\dots, x_6]$ consider the ideals $I_1=(x_1, x_5, x_6)\cap(x_4, x_5, x_6)=(x_1x_4, x_5, x_6)$ and $I_2=(x_1, x_2, x_3)$ and set $I=I_1\cap I_2=(x_1x_4, x_1x_5, x_1x_6, x_2x_5, x_2x_6, x_3x_5, x_3x_6)$. Then Corollary \ref{corollary} applies to $I_1\cap I_2$ with $u=v=3$ and $s=1$, so that $I$ is generated, up to radical, by 
\begin{eqnarray*}
\gamma_1&=&x_1x_5,\\
\gamma_2&=&x_2x_5+x_1x_6,\\
 \gamma_3&=&x_1x_4+x_3x_5+x_2x_6,\\
\gamma_4&=&x_3x_6.
\end{eqnarray*}
\noindent
Thus ara\,$I\leq 4$. We show that equality holds: to this end we exploit the inequality 
\begin{equation}\label{cd}{\rm cd}\,I\leq\,{\rm ara}\,I\end{equation}
\noindent
 (see \cite{Hu}, Theorem 3.4, or \cite{H}, Example 2, pp.~414--415), where cd denotes the so-called local cohomological dimension of $I$, which is defined as follows:
$${\rm cd}\,I=\max\{i\vert H^i_I(R)\ne 0\}.$$
\noindent
We prove that $H^4_I(R)\ne0$. This will imply that $4\leq{\rm cd}\,I\leq\,{\rm ara}\,I\leq 4$, so that equality holds everywhere. In particular, $I$ is not a s.t.c.i.. We have the following Mayer-Vietoris sequence (see \cite{Hu}, p.~15):
\begin{equation}\label{MV} \cdots \to H^4_{I_1}(R)\oplus H^4_{I_2}(R)\to H^4_{I_1\cap I_2}(R)\to H^5_{I_1+I_2}(R)\to H^5_{I_1}(R)\oplus H^5_{I_2}(R)\to\cdots\end{equation}
\noindent
In view of (\ref{cd}), since $I_1$ and $I_2$ are generated by three elements, we have that $H^i_{I_1}(R)=H^i_{I_2}(R)=0$ for all $i>3$. Moreover, $I_1+I_2=(x_1, x_2, x_3, x_5, x_6)$ is generated by a regular sequence of length 5, so that by \cite{H}, Example 2, pp.~414--415, cd\,$I_1+I_2=5$ and, on the other hand, by \cite{Hu}, Proposition 2.8, $H^i_{I_1+I_2}(R)=0$ for all $i<5$.  Thus $H^i_{I_1+I_2}(R)\ne0$ if and only if $i=5$. From (\ref{MV}) it thus follows that $H^4_{I_1\cap I_2}(R)\ne0$, as required. 
}\end{example} 
\begin{example}\label{example3}{\rm In any commutative ring $R$ with identity consider the ideal
$$I=(\mu_1, \mu_2, \mu_3)\cap(\nu_1, \nu_2, \nu_3)\cap(\xi_1, \xi_2, \xi_3),$$
\noindent
where we assume that $\mu_1\in(\nu_1)$ and $\mu_2\in(\xi_1)$. Corollary \ref{corollary} applies to the  ideal $I'=(\mu_1, \mu_2, \mu_3)\cap(\nu_1, \nu_2, \nu_3)$ with $u=v=3$ and $s=1$. Hence, according to (\ref{statement1}), it is generated, up to radical, by 4 elements $\gamma_1, \gamma_2, \gamma_3, \gamma_4\in R$ such that $\gamma_1\in (\mu_2)$; consequently, $\gamma_1\in(\xi_1)$. Thus $I''=(\gamma_1,\gamma_2, \gamma_3, \gamma_4)\cap(\xi_1, \xi_2, \xi_3)$ fulfills the assumption of Corollary \ref{corollary} for $u=4$, $v=3$ and $s=1$. Note that $\sqrt{I''}=\sqrt{I}$. Therefore ara\,$I\leq 5$. For instance, in the polynomial ring $K[x_1,\dots, x_{10}]$, where $K$ is a field, consider the ideal
$$I=(x_1x_4, x_2x_5, x_3x_6)\cap(x_1, x_7, x_8)\cap(x_2, x_9, x_{10}).$$
\noindent
It is generated, up to radical, by the following 5  polynomials:
\begin{eqnarray*}
&&(x_2x_5x_7+x_1x_3x_6)x_2,\\
&&(x_2x_5x_7+x_1x_3x_6)x_9+(x_1x_4+x_2x_5x_8+x_3x_6x_7)x_2,\\
&&(x_2x_5x_7+x_1x_3x_6)x_{10}+(x_1x_4+x_2x_5x_8+x_3x_6x_7)x_9+x_2x_3x_6x_8,\\
&&x_1x_2x_5+(x_1x_4+x_2x_5x_8+x_3x_6x_7)x_{10}+x_3x_6x_8x_9,\\
&&x_3x_6x_8x_{10}.
\end{eqnarray*}
As was proven in \cite{L3}, the cohomological dimension of a monomial ideal is equal to its projective dimension (pd). In our case, a computation with CoCoA \cite{C} yields that pd\,$I=5$ if char\,$K=0$. Hence, in view of (\ref{cd}), we have that ara\,$I=5$ in characteristic zero. Another computation with CoCoA yields that the minimum number of generators of $I$ is 15. 
}
\end{example}
\section{An application}
In this section we apply the results of the previous section to the computation of  the arithmetical ranks of a certain class of ideals generated by monomials in a polynomial ring over a field $K$. Since the arithmetical rank is the same up to radical, without loss of generality, we can restrict our attention to ideals generated by {\it squarefree} monomials. We first prove one general result. 
\begin{proposition}\label{Proposition}
Let $h,t$ be integers, $h\leq t$, and let $a_0,a_1,\dots, a_h$ be integers with $0=a_0<a_1<\cdots<a_h\leq t$. Let $x_1,\dots, x_t, y_1,\dots, y_{a_h}\in R$ and consider the following ideals of $R$:
$$I_0=(x_1,\dots, x_t),$$
\noindent and, for all indices $i$ with $1\leq i\leq h$, 
$$I_i=(y_1,\dots, y_{a_i}, x_{a_i+1},\dots, x_t).$$
\noindent
Then there are $\gamma_1,\dots, \gamma_{t+a_h-h}\in R$ generating $I_0\cap I_1\cap\cdots\cap I_h$ up to radical such that
\begin{equation}\label{equation1} [\gamma_1,\dots, \gamma_{a_h}]\subset [y_1,\dots, y_{a_h}],\end{equation}
\noindent and
$$\gamma_i=x_{i-a_h+h}$$
\noindent
for all indices $i$ with $2a_h-h+1\leq i\leq t+a_h-h$.
\end{proposition} 
\demo We proceed by induction on $h$. The claim is trivially true for $h=0$. Suppose that $h\geq 1$ and that the claim is true for all smaller indices. There are $\bar\gamma_1,\dots, \bar\gamma_{t+a_{h-1}-h+1}\in R$ generating $I_0\cap I_1\cap \cdots\cap I_{h-1}$ up to radical such that
\begin{equation}\label{13}[\gamma_1,\dots, \gamma_{a_{h-1}}]\subset [y_1,\dots, y_{a_{h-1}}],\end{equation}
\noindent and
$$\gamma_i=x_{i-a_{h-1}+h-1}$$
\noindent
for all indices $i$ with $2a_{h-1}-h+2\leq i\leq t+a_{h-1}-h+1$. Note that 
\begin{eqnarray}\label{radical} &&\sqrt{I_0\cap\cdots\cap I_{h-1}\cap I_h}\nonumber\\
&=&\sqrt{(\bar\gamma_1,\cdots, \bar\gamma_{t+a_{h-1}-h+1})\cap(y_1,\dots, y_{a_h}, x_{a_h+1},\dots, x_t)}\nonumber\\
&=&\sqrt{(\bar\gamma_1,\dots, \bar\gamma_{a_h+a_{h-1}-h+1},x_{a_h+1},\dots, x_t)\cap(y_1,\dots, y_{a_h}, x_{a_h+1},\dots, x_t)}\nonumber\\
&&\mbox{
(note that $a_h+a_{h-1}-h+1\geq 2a_{h-1}-h+2$)}\nonumber
\\
&=&\sqrt{(\bar\gamma_1,\cdots, \bar\gamma_{a_h+a_{h-1}-h+1})\cap(y_1,\dots, y_{a_h})+(x_{a_h+1},\dots, x_t)}.
\end{eqnarray}
\noindent
In view of (\ref{13}), Corollary \ref{corollary} applies to $(\bar\gamma_1,\cdots, \bar\gamma_{a_h+a_{h-1}-h+1})\cap(y_1,\dots, y_{a_h})$ for  $u=a_h+a_{h-1}-h+1$, $v=a_h$ and $s=a_{h-1}$.  Hence, by (\ref{statement2}),  this ideal is generated up to radical by $\gamma_1,\dots, \gamma_{2a_h-h}\in R$ such that
$$[\gamma_1,\dots, \gamma_{a_h}]\subset [y_1,\dots, y_{a_h}].$$
\noindent
Set $\gamma_i=x_{i-a_h+h}$ for all indices $i$ such that $2a_h-h+1\leq i\leq t+a_h-h.$ Then, by (\ref{radical}), $\gamma_1,\dots, \gamma_{t+a_h-h}$ generate $I_0\cap \cdots\cap  I_{h-1}\cap I_h$ up to radical. This completes the proof.
\par\medskip\noindent
\begin{corollary}\label{corollary2} Let $K$ be a field and let $x_1,\dots, x_t, y_1,\dots, y_t$ be pairwise distinct indeterminates over $K$. In the polynomial ring $R=K[x_1,\dots, x_t, y_1,\dots, y_t]$ consider the ideal
$$I=(x_1,\dots, x_t)\cap(y_1,x_2,\dots, x_t)\cap\cdots\cap(y_1,\dots, y_{t-1}, x_t)\cap(y_1,\dots, y_t).$$
\noindent
Then $I$ is a s.t.c.i..
\end{corollary}
\demo It suffices to apply Proposition \ref{Proposition} for $h=t$ (which implies that $a_i=i$ for all $i=1,\dots, t$). This yields ara\,$I\leq t$. Since height\,$I=t$, equality holds. This completes the proof.
\par\medskip\noindent
\begin{remark}\label{remark2}{\rm Recall that the ideals generated by squarefree monomials in a polynomial ring over a field are the face ideals of simplicial complexes: we refer to Bruns and Herzog \cite{BH}, Section 5, for the basic notions on this topic.  The ideal $I$ in the claim of Corollary \ref{corollary2} is the face ideal of the simplicial complex $\Delta$ on the vertices $x_1,\dots, x_t, y_1,\dots, y_t$ whose facets are
\begin{eqnarray*}
F_0&=&\{y_1,\dots, y_t\},\\
F_1&=&\{x_1,y_2,\dots, y_t\},\\
\vdots&&\\
F_{i-1}&=&\{x_1,\dots,x_{i-1}, y_i,\dots, y_t\}\\
F_i&=&\{x_1,\dots, x_i, y_{i+1}, \dots, y_t\},\\
\vdots&&\\
F_t&=&\{x_1,\dots, x_t\}.
\end{eqnarray*}
\noindent
For all indices $i>0$,  
$$F_i\cap\left(\displaystyle\cup_{i=0}^{i-1}F_i\right)=\{x_1,\dots, x_{i-1}, y_{i+1},\dots, y_t\},$$
\noindent
which is a maximal proper subset of $F_{i-1}$. According to \cite{BH}, Definition 5.1.11., $\Delta$ is a so-called {\it shellable} complex. By virtue of \cite{BH}, Theorem 5.1.13., as a consequence, we have that $R/I$ is a Cohen-Macaulay ring. By the results in \cite{BM}, $I$  also defines a reducible variety of minimal degree according the classification given in \cite{EG}, Section 4.  Another class of varieties of minimal degree which are s.t.c.i. (and whose defining ideal is generated by a set formed of both monomials and binomials of degree two) has been recently characterized in \cite{B2}.  
}
\end{remark}

\end{document}